\documentclass{article}

\usepackage{amsmath}
\usepackage{amsthm}
\usepackage[utf8]{inputenc}
\usepackage{amsfonts}
\usepackage{hyperref}
\hypersetup{
	colorlinks = true,
	linkcolor = purple, 
	filecolor = purple, 
	citecolor = purple,
	urlcolor = purple
}
\usepackage{amssymb}
\usepackage{graphicx}
\usepackage{stackrel}
\usepackage[english]{babel}
\usepackage{url}
\usepackage{xcolor}
\usepackage{tikz-cd}
\usepackage{pdfpages}
\usepackage{csquotes}

\usepackage[normalem]{ulem}
\usepackage{soul}
\allowdisplaybreaks

\newtheorem{theorem}{Theorem}[section]
\newtheorem{introtheorem}{Theorem}
\newtheorem{prop}[theorem]{Proposition}
\newtheorem*{prop*}{Proposition}
\newtheorem{lem}[theorem]{Lemma}
\newtheorem{cor}[theorem]{Corollary}

\newtheorem{introconj}[introtheorem]{Conjecture}
\newtheorem*{conj*}{Conjecture}
\theoremstyle{definition}
\newtheorem{defn}[theorem]{Definition}
\newtheorem{nota}[theorem]{Notation}
\newtheorem{constr}[theorem]{Construction}

\newtheorem{rem}[theorem]{Remark}

\begin{document}
	\title{Low degree motivic Donaldson-Thomas invariants of the three-dimensional projective space}
	\author{Anna M. Viergever}
	\date{ }
	\maketitle 
	
\begin{abstract}
	Levine has constructed motivic analogues of virtual fundamental classes, living in cohomology of Witt sheaves. We use this to define motivic Donaldson-Thomas invariants $\tilde{I}_n$ for $\mathbb{P}^3$ over $\mathbb{R}$. We show that for $n$ odd, $\tilde{I}_n = 0$ and we compute $\tilde{I}_2 = 10, \tilde{I}_4 = 25$ and $\tilde{I}_6 = -50$. We then make a conjecture about the general case, which could be a motivic analogue of a classical theorem of Maulik-Nekrasov-Okounkov-Pandharipande. 
\end{abstract}
	
	\tableofcontents 
	
	\phantomsection
	\addcontentsline{toc}{section}{Introduction}
	
	\section*{Introduction} 
	\subsection*{Motivation and the theorem of Maulik-Nekrasov-Okounkov-Pandharipande}
	If one tries to count a certain type of objects for which there exists a reasonable moduli space, one could do that by doing intersection theory on that moduli space. However, these moduli spaces often have all sorts of bad singularities, or may not have the expected dimension. Virtual fundamental classes, introduced by Behrend and Fantechi in their paper \cite{BehrendFantechiISNC}, provide a way to make reasonable computations in spite of these problems. For example, if we try to count ideal sheaves on a smooth projective scheme which are of a given length, the corresponding moduli space is the Hilbert scheme. The degrees of the corresponding virtual fundamental classes in this situation are called \textit{Donaldson-Thomas invariants}. They were first constructed and defined by Donaldson and Thomas in \cite{DonaldsonThomas} and \cite{ThomasHCI}.	
	
	A particular example of a computation of Donaldson-Thomas invariants is done in the papers \cite{MaulikI} and \cite{MaulikII}, by Maulik, Nekrasov, Pandharipande and Okounkov. If we take a smooth projective toric threefold $X$ over $\mathbb{C}$ with an action of the three dimensional torus $\mathbb{T}$ on it, we can look at the Hilbert scheme $\text{Hilb}^n(X)$ of ideal sheaves on $X$ of length $n$. To this, one associates a virtual fundamental class, of which the degree $I_n\in\mathbb{Z}$ is the Donaldson-Thomas invariant.	In \cite{MaulikI} and \cite{MaulikII}, there is a proof of the fact that 
	\begin{equation} \label{MNOP result}
		\sum_{n\geq 0}I_nq^n = M(-q)^{\deg(c_3(T_X\otimes K_X))}
	\end{equation} where $M(q) = \prod_{n\geq 1}(1-q^n)^{-n}$ is the MacMahon function (see MacMahon's book \cite[Article 43]{MacMahon} and Stanley's book \cite[Corollary 7.20.3]{Stanley}), $T_X$ is the tangent bundle on $X$ and $K_X$ is the canonical line bundle of $X$. 
	
	The proof uses the virtual localization formula by Graber and Pandharipande, see \cite{GraberLVC}, to compute the degree of the $\mathbb{T}$-equivariant virtual fundamental class instead. This degree  lives in $\text{CH}^*(B\mathbb{T}) \cong \mathbb{Z}[x,y,z]$ where, $B\mathbb{T}$ is the classifying space of $\mathbb{T}$, see Totaro's paper \cite{TotaroCRCS}. As the degree of the equivariant virtual fundamental class lands $\text{CH}^0(B\mathbb{T})\cong\mathbb{Z}$, one can deduce the Donaldson-Thomas invariant $I_n$ from it. In this case, the virtual localization formula gives that
	$$I_n = \sum_{[\mathcal{I}]\in \text{Hilb}^n(X)^{\mathbb{T}}} \frac{e(\text{Ext}^2(\mathcal{I},\mathcal{I}))}{e(\text{Ext}^1(\mathcal{I},\mathcal{I}))}.$$
	Here, the sum is over all ideal sheaves $\mathcal{I}$ on $X$ of length $n$ which are fixed under the induced action of $\mathbb{T}$ on $\text{Hilb}^n(X)$. Such ideals are locally given by monomial ideals. The authors then compute the trace, as a $\mathbb{T}$-representation, of the virtual tangent space $$\text{Ext}^1(\mathcal{I},\mathcal{I}) - \text{Ext}^2(\mathcal{I},\mathcal{I})$$ of each fixed ideal sheaf $\mathcal{I}$. From this trace, one can read off the corresponding equivariant virtual Euler class and express this in terms of the Euler classes of the standard line bundles $\mathcal{O}(1,0,0), \mathcal{O}(0,1,0)$ and $\mathcal{O}(0,0,1)$ on $B\mathbb{T}$. To deduce the above formula, the Bott residue formula \cite{Bott} is used. 
	
	Maulik, Nekrasov, Pandharipande and Okounkov conjectured that the formula (\ref{MNOP result}) would hold for the Donaldson-Thomas invariants of all smooth projective threefolds, not just the toric ones. This conjecture was proven by Li, see \cite{Li}, and by Levine-Pandharipande, see \cite{LevinePandharipande}, using completely different methods.
	
	\subsection*{Main result}
	We study an analogue of these results for $\mathbb{P}^3$, replacing the integer-valued Donaldson-Thomas invariants $I_n$ with a parallel version constructed using motivic virtual fundamental classes with values in cohomology of Witt sheaves defined by Levine in \cite{LevineVLEWC}, and replacing the Graber-Pandharipande virtual localization with the corresponding virtual localization formula proven by Levine in the same paper. To use this localization formula, one has to invert an element, which is in degree $d\geq 1$ of the cohomology groups of Witt sheaves on the classifying space of the torus used. 
	
	The quadratic degrees of these motivic Donaldson-Thomas invariants will land in the Witt ring of the base field. As $\mathcal{W}(\mathbb{C}) \cong \mathbb{Z}/2\mathbb{Z}$, one would only recover the classical value modulo two in the case of the complex numbers. For $\mathbb{F}$ a finite field of odd characteristic, one has $4\mathcal{W}(\mathbb{F}) = 0$ which means that the computation yields no information if $2$ has to be inverted. On the other hand, $\mathcal{W}(\mathbb{R}) = \mathbb{Z}$ via the signature map, so working over $\mathbb{R}$ will give an invariant that one can view as a real version of the classical Donaldson-Thomas invariants; we therefore take the base field to be $\mathbb{R}$ although most of our arguments will work for any field which has an embedding into $\mathbb{R}$.
	
	Having an action by the torus $\mathbb{T}$ is not useful for Witt-sheaf localization. Indeed, since 
	$$H^d(B\mathbb{T},\mathcal{W}) = \begin{cases}
		\mathcal{W}(k) &\text{ if } d=0\\
		0 &\text{otherwise}
	\end{cases} $$
	$\mathbb{T}$-localization would require inverting $0\in H^{*\ge1}(B\mathbb{T},\mathcal{W})$, yielding no useful information. Following the discussion in \cite{LevineABLEWC}, using  the normalizer $N_S$ of the torus in $\text{SL}_2$ yields a useful version of localization applicable to Witt-sheaf cohomology.
	
	We recall that $N_S$ is the subgroup of $\text{SL}_2$ generated by matrices
	$$\left(\begin{matrix} t & 0 \\ 0 & t^{-1}\end{matrix} \right) \text{ and } \sigma = \left(\begin{matrix} 0 & 1 \\ -1 & 0\end{matrix} \right).$$
	Let $a,b\in\mathbb{Z}$ be suitably general odd integers. There is an action of $N_S$ on $\mathbb{P}^3$ given by 
	\begin{align*}
		\left(\begin{matrix} 
			t & 0\\
			0 & t^{-1}
		\end{matrix} \right) &\cdot [X_0:X_1:X_2:X_3] = [t^aX_0: t^{-a}X_1: t^bX_2: t^{-b}X_3]\\
		\sigma &\cdot [X_0:X_1:X_2:X_3] = [-X_1:X_0:-X_3:X_2].
	\end{align*}
	This action does not have fixed points on $\mathbb{P}^3$, but there are the two fixed couples $\{[1:0:0:0], [0:1:0:0] \}$ and $\{[0:0:1:0], [0:0:0:1] \}$. 
	
	Levine \cite{LevineOST} has constructed an orientation for $\text{Hilb}^n(X)$ for $X$ any smooth projective threefold with a given isomorphism $K_X\cong L^{\otimes 2}$ for some invertible sheaf $L$ on~$X$. This implies that there is a well defined quadratic degree $\tilde{I}_n\in \mathcal{W}(\text{Spec}(\mathbb{R}))\cong\mathbb{Z}$ of the virtual fundamental class of $\text{Hilb}^n(\mathbb{P}^3)$.
	
	Using \cite[Theorem 6.7]{LevineVLEWC}, if the torus inside $N_S$ acts with isolated fixed points, we have that
	\begin{align*}
		\tilde{I}_n &= \sum_{[\mathcal{I}]\in \text{Hilb}^n(\mathbb{P}^3)^{N_S}} \frac{\widetilde{\deg}(e(\text{Ext}^2(\mathcal{I},\mathcal{I})))}{\widetilde{\deg}(e(\text{Ext}^1(\mathcal{I},\mathcal{I})))}.
	\end{align*}
	where $e(-)$ is the $N_S$-equivariant Euler class and $\widetilde{\deg}$ denotes the $\mathcal{W}(\mathbb{R})$-valued quadratic degree map, and we identify $\mathcal{W}(\mathbb{R})$ with $\mathbb{Z}$ via the signature. This already tells us that our quadratic invariant $\tilde{I}_n$ will vanish for odd $n$, and  we can concentrate on the case of even $n$.	
	
	We first compute the trace of the virtual representation $\text{Ext}^2(\mathcal{I},\mathcal{I}) - \text{Ext}^1(\mathcal{I},\mathcal{I})$ of $N_S$ for ideal sheaves $\mathcal{I}$ that are fixed by the $N_S$-action, using the strategy of \cite{MaulikI} and \cite{MaulikII}. We can compute the Euler classes from that in terms of the Euler classes of certain canonical rank two bundles using the results \cite[Proposition 5.5 and Theorem 7.1]{LevineMEC}. We then apply this to compute the first three nonzero real Donaldson-Thomas invariants. 
	\begin{introtheorem}\label{main theorem}
		Let $X = \mathbb{P}^3_{\mathbb{R}}$. For $n\geq 0$, let $\tilde{I}_n$ be the quadratic degree of the motivic virtual fundamental class associated to~$\text{Hilb}^n(X)$. Then $\tilde{I}_2 = 10, \tilde{I}_4 = 25$ and $\tilde{I}_6 = -50$. 
	\end{introtheorem} 
	For $n=8$ and higher, the present method does not work to compute the corresponding motivic Donaldson-Thomas invariant. Namely, using coordinates $x=X_1/X_0,y=X_2/X_0,z=X_3/X_0$ on the open $U_0$ defined by $\{X_0\neq 0\}\subset\mathbb{P}^3$,  the  length eight  $\sigma$-invariant ideal sheaves whose restriction to $U_0$ are given by $(x+\lambda yz,y^2,z^2)$ for $\lambda\in \mathbb{R}^*$ are all fixed by $N_S$, so the $N_S$-action on $\text{Hilb}^8(X)$ has non-isolated fixed points. This propagates to give a positive dimensional $N_S$-fixed locus on $\text{Hilb}^n(X)$ for all even $n\ge8$. Although the localization method does in principle work even in the case of non-isolated fixed points, the computations become more difficult, and we will not consider these cases in this paper.	
	
	The following expression involving the MacMahon function
	\begin{equation}\label{eqn:MacMahon}
		M(-q^2)^{-10} = 1 + 10q^2 + 25q^4 - 50q^6 - 240q^8 + \cdots 
	\end{equation}
	has coefficients of $q^{2n}$ agreeing with $\tilde{I}_{2n}$ for our computed values $n=1,2,3$. 
	Note that, for $X$ a smooth projective threefold over a field $k$ of characteristic different from 2, the obstruction sheaf $\text{ob}_2(X)$ on $\text{Hilb}^2(X)$ is a locally free sheaf whose fiber at $[\mathcal{I}]\in\text{Hilb}^2(X)$ is  the ``traceless'' summand $\text{Ext}^2_0(\mathcal{I}, \mathcal{I})$ of $\text{Ext}^2(\mathcal{I}, \mathcal{I})$. The same heuristic argument used in \cite{MaulikII} to justify the exponent $c_3(T_X\otimes K_X)$ in the formula (\ref{MNOP result}) leads one to expect the exponent $\widetilde{\deg}(e(\text{ob}_2(X)))$ in a quadratic replacement. Together with \eqref{eqn:MacMahon} and the fact that $\widetilde{\deg}(e(\text{ob}_2(\mathbb{P}^3)))=-10$, this leads us to the following conjecture.
	\begin{introconj}
		Let $X$ be a smooth projective threefold over $\mathbb{R}$, equipped with an isomorphism  $K_X\cong L^{\otimes 2}$ for $L$ some line bundle on $X$. For $n\geq 0$, let $\tilde{I}_n$ be the  quadratic degree of the motivic virtual fundamental class associated to~$\text{Hilb}^n(X)$. Then 
		$$\sum_{n\geq 0}\tilde{I}_nq^n = M(-q^2)^{-\widetilde{\deg}(e(\text{ob}_2(X)))}.$$
	\end{introconj} 
	
	We would like to highlight that Espreafico and Walcher \cite{EspreaficoDT} computed the full generating series of a version of the motivic Donaldson-Thomas invariants of $\mathbb{A}^3$ over a characteristic zero field. Their strategy is to first consider the motivic generating series for the Hilbert schemes of $\mathbb{A}^3$, which was computed by Behrend, Bryan and Szendr\"oi in \cite{BehrendMDT}. Here ``motivic'' means with values in $K_0(\text{Var}_k)$,  the Grothendieck ring of varieties over a base field $k$.  They then apply a version of the compactly supported $\mathbb{A}^1$-Euler characteristic, which was constructed by Arcila-Maya, Bethea, Opie, Wickelgren and Zakharevich in \cite{WickelgrenEulercharacteristic}, as a ring homomorphism $\chi_c: K_0(\text{Var}_k)\to \text{GW}(k)$, where  $\text{GW}(k)$ denotes the Grothendieck-Witt ring of $k$. This gives a series which refines known results over $\mathbb{C}$ and $\mathbb{R}$. 
	
	One cannot compute the motivic Donaldson-Thomas invariants for $\mathbb{P}^3$ using this method directly. The reason is that for $X$ a Calabi-Yau manifold and $n\geq 1$, the Donaldson-Thomas invariant $I_n$ is equal to the Behrend weighted Euler characteristic of $\text{Hilb}^n(X)$ up to a sign, which in turn is equal to the usual Euler characteristic (again up to sign) in case $\text{Hilb}^n(X)$ is smooth; see for instance \cite[Theorem 4.12]{Behrend-FantechiSOTHSPT}. However, this no longer holds for threefolds which are not Calabi-Yau. For instance, let $\mathbb{L}$ be the class of $\mathbb{A}^1$ in the Grothendieck group of varieties over $k$, then the motivic (in the sense of $K_0(\text{Var}_k)$) generating series of Donaldson-Thomas invariants for $\mathbb{P}^3$ of \cite{BehrendMDT} is 
	$$1 - (\mathbb{L}^{-\frac{3}{2}} + \mathbb{L}^{-\frac{1}{2}} + \mathbb{L}^{\frac{1}{2}} + \mathbb{L}^{\frac{3}{2}})t + (\mathbb{L}^{-3} + 2\mathbb{L}^{-2} + 4\mathbb{L}^{-1} + 4 + 4\mathbb{L} + 2\mathbb{L}^2 + \mathbb{L}^3)t^2 + \cdots$$ 
	in $K_0(\text{Var}_k)[\mathbb{L}^{-\frac{1}{2}}]$. Evaluating the series with Behrend's weighted motivic Euler characteristic maps $\mathbb{L}^{-\frac{1}{2}}$ to $-1$, meaning that the first coefficient is mapped to $-4$ and the second coefficient is mapped to $18$; as expected, these agree (up to sign) with the topological Euler characteristics of $\mathbb{P}^3$ and of $\text{Hilb}^2(\mathbb{P}^3)$, respectively.  However, following the formula (\ref{MNOP result}), the generating series for the classical Donaldson-Thomas invariants for $\mathbb{P}^3_{\mathbb C}$ is 
	\[
	M(-q)^{-20}=1+20q+150q^2+\ldots
	\]
	so $I_1(\mathbb{P}^3_{\mathbb{C}})=20$, $I_2(\mathbb{P}^3_{\mathbb{C}})=150$. Similarly, applying $\chi_c$ to   the $K_0(\text{Var}_{\mathbb{R}})$-series, one finds the form $10\langle -1 \rangle + 8\langle 1 \rangle$ for the second coefficient, which does not have signature $-10 = \tilde{I}_2$. Thus, as it stands, the approach of Espreafico and Walcher needs at least to be suitably modified to yield formulas for the quadratic, or even real, Donaldson-Thomas invariants for non-Calabi-Yau threefolds.
	
	The results presented here also form a chapter in the authors thesis, which was
	submitted on May 30’th, 2023.
	
	For the computations of $\tilde{I}_2, \tilde{I}_4$ and $\tilde{I}_6$, some SAGE-code has been used, which can be retrieved on the author's webpage, see \href{https://sites.google.com/view/anneloes-viergever/research}{https://sites.google.com/view/anneloes-viergever/research}.
	
	\subsection*{Acknowledgements} 
	I deeply thank my advisor Marc Levine for suggesting this project to me and for all of his invaluable help over the last years. 
	
	I would also like to thank Alessandro d'Angelo, for teaching me almost everything I know about virtual fundamental classes, and Herman Rohrbach, for explaining me how classifying spaces work. I also wish to thank Ludvig Modin, for organizing a seminar on deformation theory with me which has been very helpful for my understanding of Hilbert schemes, and Nicolas Dupr\'e, for helpful discussions on semi-simple actions of algebraic groups.  
	
	This work was funded by the RTG Graduiertenkolleg 2553.
	
	\section{Motivic Donaldson-Thomas invariants of a smooth projective threefold}
	Let $\mathcal{E}$ be an SL-oriented motivic ring spectrum, see for instance Ananyevskiy's paper \cite{AnanyevskiySLPBT} for a definition and more details. In Levine's paper \cite{LevineISNC}, it is discussed how one can define a motivic analogue of a virtual fundamental class for a perfect obstruction theory $E_\bullet$ on a quasi-projective scheme $Z$ over a perfect field $k$. This class is an element of the Borel-Moore homology $\mathcal{E}^{B.M.}(Z,\mathbb{V}(E_\bullet))$ where $\mathbb{V}(E_\bullet)=\text{Spec}(\text{Sym}^*E_0)-\text{Spec}(\text{Sym}^*E_1)$ is the virtual vector bundle associated to $E_\bullet$. 
	
	If $E_\bullet$ has virtual rank zero and there is an isomorphism from the determinant of the obstruction theory to the square of a line bundle, a choice of such an isomorphism (called an \textit{orientation} of $E_\bullet$) defines a degree of this class as an element of $\mathcal{E}^{0,0}(k)$ (see \cite[Section 8.1]{LevineISNC}). In Levine's paper \cite{LevineVLEWC}, this definition is extended to an equivariant setting if the scheme has an action on it of a smooth closed subgroup of $\text{GL}_n$, and \cite{LevineVLEWC} provides a proof of a virtual localization formula for this situation. 
	
	In this section, we will take $\mathcal{E}$ to be cohomology of Witt sheaves. For $X$ a smooth projective threefold over a field and $n\in\mathbb{Z}_{\geq 1}$, we study the conditions under which the motivic virtual fundamental class of $\text{Hilb}^n(X)$ has a well defined quadratic degree. We will see that there is a perfect obstruction theory on $\text{Hilb}^n(X)$, and that there is an orientation if there exists an isomorphism $K_X\cong L^{\otimes 2}$ for $L$ a line bundle on $X$. Thus, in case $K_X\cong L^{\otimes 2}$, we may apply the $W(k)$-valued quadratic degree map $\widetilde{\deg}$ to the motivic virtual fundamental class for $\text{Hilb}^n(X)$, arriving at the quadratic Donaldson-Thomas invariant $\tilde{I}_n(X)\in W(k)$. We will mostly be interested in the case $k=\mathbb{R}$, giving  $\tilde{I}_n(X)\in \mathbb{Z}$ by taking the signature.
	\begin{nota} 
		Let $N_S$ be the normalizer of the diagonal torus in $\text{SL}_2$. $N_S$ is the subgroup of  $\text{SL}_2$  generated by the diagonal matrices
		$$\Bigl\{ \left(\begin{matrix} 
			t & 0\\
			0 & t^{-1}\end{matrix}\right)\Bigr\} \subset \text{SL}_2,
		$$
		\text{ and the additional matrix } $$\sigma :=  \left(\begin{matrix} 
			0 & 1\\
			-1 & 0\end{matrix}\right).$$
	\end{nota} 
	
	We also study smooth projective threefolds with an action of $N_S$, and conditions under which one can define an equivariant motivic Donaldson-Thomas invariant. 
	
	In Section \ref{section Witt cohomology}, we recall some results of Levine about the cohomology of Witt sheaves on $BN_S$ together with some computations of Euler classes of rank two bundles on $BN_S$ which will be needed later on. We then discuss the perfect obstruction theory in Section \ref{section POT} and the orientation in Section \ref{section orientation}. 
	
	\subsection{Witt cohomology of $BN_S$ and Euler classes of canonical bundles} \label{section Witt cohomology}
	For algebraic groups $G$ which are subgroups of $\text{GL}_n$, the corresponding classifying spaces $BG$ have been constructed by Totaro in \cite[Remark 1.4]{TotaroCRCS} and developed further by Edidin and Graham in \cite{EdidinGrahamEIT}. We use the same model as \cite{LevineMEC}, which comes from the model used in the paper \cite[Section 4.2.2]{MorelVoevodskyAHTS} of Morel and Voevodsky. 
	\begin{nota} 
		Let $\rho^-$ be the representation $N_S\to\mathbb{G}_m$ which sends $\sigma$ to $-1$ and the diagonal matrices to $1$.
	\end{nota} 
	\begin{nota} 
		Let $\gamma_S$ be the pullback of $\mathcal{O}(1)$ under the map $BN_S\to B\mathbb{G}_m\cong \mathbb{P}^\infty$ induced by $\rho^-$. 
	\end{nota}  
	\begin{rem}\label{remark gamma bundle}
		The Picard group of $BN_S$ is isomorphic to $\mathbb{Z}/2\mathbb{Z}$, and it is generated by $\gamma_S$. Namely, a character of the torus is of the form 
		$$\left(\begin{matrix}
			t & 0 \\ 0 & t^{-1} \end{matrix}\right)\mapsto t^{a-b}$$
		for some $a,b\in\mathbb{Z}$ and from the relation 
		\begin{equation}\label{relation sigma and t}
			\left(\begin{matrix}
				t & 0 \\ 0 & t^{-1} \end{matrix}\right)\sigma = \sigma\left(\begin{matrix}
				t^{-1} & 0 \\ 0 & t \end{matrix}\right) 
		\end{equation} 
		we see that for a character $\rho:N_S\to\mathbb{G}_m$ we have that $t^{a-b}\rho(\sigma) = \rho(\sigma)t^{b-a}$ so that $a=b$. As $\sigma^2 = -\text{Id}$ we have that $\rho(\sigma^2) = (-1)^0=1$ and so a character is either the trivial character or $\rho^-$. 
	\end{rem}
	Morel has extended the Witt ring of a field to a sheaf $\mathcal{W}$ on the category of smooth projective schemes over a given base field, see Morel's paper \cite[Chapter 2]{MorelIAHT}. This gives rise to an $\text{SL}$-oriented cohomology theory $H^*(-, \mathcal{W})$. For more details see \cite[Section 2]{LevineAEGQF}. 
	
	The Witt sheaf cohomology of $BN_S$ has been computed by Levine.
	\begin{nota}
		Let $\bar{q}\in H^0(BN_S,\mathcal{W})$ be the element constructed in \cite[Section 5]{LevineMEC}. Let $p:BN_S\to B\text{SL}_2$ be the map constructed in \cite[Section 5]{LevineMEC}. Let $\mathcal{T}$ denote the tangent bundle of $BN_S$ over $B\text{SL}_2$. 	Let $e$ be the Euler class of the canonical rank two bundle $\mathbb{A}^2\times^{\text{SL}_2}B\text{SL}_2$ on $B\text{SL}_2$, where $\text{SL}_2$ acts on $\mathbb{A}^2$ by right matrix multiplication. 
	\end{nota} 
	\begin{prop}[\cite{LevineMEC}, Proposition 5.5]\label{Marcs computation of cohomology of W}
		Let $k$ be a perfect field. Let $W(k)[x_0,x_2]$ be the graded polynomial algebra over $W(k)$ on the generators $x_0$ of degree zero and $x_2$ of degree~$2$. Then $x_0\mapsto \langle \bar{q}\rangle, x_2\mapsto p^*e$ defines a $W(k)$-algebra isomorphism 
		$$\psi: W(k)[x_0,x_2]/(x_0^2-1,(1+x_0)x_2)\to H^*(BN_S,\mathcal{W}).$$
		Moreover, $H^{*\geq 2}(BN_S,\mathcal{W}(\gamma_S))$ is the quotient of the free $H^{*}(BN_S,\mathcal{W})$-module on the generator $e(\mathcal{T})$ modulo the relation $(1+ \langle q\rangle)e(\mathcal{T}) = 0$.  
	\end{prop}
	In \cite{LevineMEC}, there is also a computation of the Euler classes of all rank two vector bundles on $BN_S$. Following \cite[Section 6]{LevineMEC}, for $m\geq 1$, consider the representation $\rho_m: N_S\to \text{GL}_2(\mathbb{A}^2)$ given by 
	\begin{align*}
		\left(\begin{matrix} 
			t & 0\\
			0 & t^{-1}
		\end{matrix} \right) &\mapsto\left(\begin{matrix} 
			t^m & 0\\
			0 & t^{-m}
		\end{matrix} \right) \\
		\sigma &\mapsto \left(\begin{matrix} 
			0 & 1\\
			(-1)^m & 0
		\end{matrix} \right)
	\end{align*}
	and let $\rho_m^-$ be given by $\rho_m^-(\sigma) = - \rho_m(\sigma)$. Finally, let $\rho_0$ be the trivial representation and let $\rho_0^- = \rho^-$. The representations $\rho_m$ give rise to rank two vector bundles $\tilde{O}(m)$ on $BN_S$ and the representations $\rho_m^-$ similarly give rise to rank two vector bundles $\tilde{O}^-(m)$ on $BN_S$; we sometimes write $\rho_m^+$ for $\rho_m$. We orient $\rho_m^\pm$ by taking $e_1\wedge e_2$ as the chosen generator of $\det\rho_m^\pm$, defining a canonical isomorphism $\det\rho_m^\pm\cong \rho_0$ for $m$ odd and $\det\rho_m^\pm\cong \rho_0^-$ for $m$ even.

	By \cite[Lemma 6.1]{LevineMEC}, $\tilde{O}(2)$ is isomorphic to $\mathcal{T}$. Levine has proven the following result. 
	\begin{theorem}[\cite{LevineMEC}, Theorem 7.1]\label{Marcs computation of euler classes}
		Suppose that $k$ is a field of characteristic zero or of characteristic $p>2$ such that $p$ and $m$ are coprime. Then 
		$$e(\tilde{O}(m)) = \begin{cases}
			m\cdot p^*e \in H^2(BN_S,\mathcal{W}) &\text{ if } m\equiv 1 \text{ mod } 4\\
			-m\cdot p^*e \in H^2(BN_S,\mathcal{W}) &\text{ if } m\equiv 3 \text{ mod } 4\\
			\frac{m}{2}\cdot e(\mathcal{T}) \in H^2(BN_S,\mathcal{W}(\gamma_S)) &\text{ if } m\equiv 2 \text{ mod } 4\\
			-\frac{m}{2}\cdot e(\mathcal{T}) \in H^2(BN_S,\mathcal{W}(\gamma_S)) &\text{ if } m\equiv 0 \text{ mod } 4
		\end{cases}$$
		Furthermore, $e(\tilde{O}^-(m)) = -e(\tilde{O}(m))$. 
	\end{theorem}
	\begin{rem}\label{remark NS representations are semisimple}
		Note that the above theorem gives a complete computation of Euler classes of bundles on $BN_S$. We have considered all possible twists as we saw in Remark \ref{remark gamma bundle}. Also, let $\rho:N_S\to \text{GL}(V)$ be a representation of $N_S$. Then we can restrict this to a representation of the torus inside of $\text{SL}_2$, and as representations of the torus $T$ are semisimple, this will be diagonalizable. Adding the element $\sigma$ and using the relation (\ref{relation sigma and t}), we see that $\rho$ is a direct sum of representations $\rho_m$ for $m\geq 0$ or $\rho^-$ (noting that any trivial subrepresentation is a $\rho_0$ and that $\rho_m$ and $\rho_m^-$ are isomorphic as representations). Note that this also proves that $N_S$-representations are semi-simple. 
	\end{rem}
	\begin{rem} \label{remark sign Euler classes}
		The Euler class $e(\tilde{O}^\pm(m))$ depends on a choice of isomorphism $\det(e(\tilde{O}^\pm(m)))\to \mathcal{O}_{BN_S}$ if $m$ is odd, or $\det(e(\tilde{O}^\pm(m)))\to \gamma_S$ if $m$ is even. As the isomorphism of representations $\tilde{O}^-(m)\to \tilde{O}(m)$ has determinant $-1$, one finds the minus sign in the formula for $e(\tilde{O}^-(m))$.  
	\end{rem}
	
	\subsection{Perfect obstruction theory}\label{section POT}
	We recall the following definition. 
	\begin{defn} 
		Let $Y$ be a scheme over a field $k$, and let $\mathbb{L}_{Y}$ be the cotangent complex of $Y$ over $k$. A \textit{perfect obstruction theory} on $Y$ is a perfect complex $E_\bullet$ supported in cohomological degrees $0$ and $-1$, together with a morphism $E_\bullet\to \mathbb{L}_{Y}$, which defines an isomorphism on cohomology groups in degree $0$ and a surjection on cohomology groups in degree $-1$. 
	\end{defn} 
	Let $X$ be a smooth projective threefold over a field $k$. Let $n \geq 1$ and let $\text{Hilb}^n(X)$ be the Hilbert scheme of ideal sheaves of length $n$ on $X$. 
	
	Thomas has constructed a perfect obstruction theory on $\text{Hilb}^n(X)$ for $X$ a smooth projective threefold. This construction of a perfect obstruction theory on $\text{Hilb}^n(X)$  for smooth projective threefolds has been rephrased in somewhat more invariant terms by Behrend-Fantechi in \cite[Section 1.5]{Behrend-FantechiSOTHSPT}.
	\begin{rem} 	
		Let $E_\bullet$ denote the perfect obstruction theory described above. If $X$ has a linear $N_S$-action and the very ample invertible sheaf $\mathcal{O}_X(1)$ defining the projective embedding has a $N_S$-linearization, one can show that the action extends to an action of $N_S$ on $\text{Hilb}^n(X)$ and this gives an $N_S$-linearization on $E_\bullet$. The perfect obstruction theory $E_\bullet$ now gives rise to an $N_S$-equivariant perfect obstruction theory. See \cite[Theorem 6.4]{LevineOST} or \cite[Proposition 2.4]{Behrend-FantechiSOTHSPT} for details.   
	\end{rem}
	\begin{rem} 
		Suppose that $H^i(X,\mathcal{O}_X) =0$ for $i\geq 1$. Let $E_\bullet$ denote the perfect obstruction theory described above. Because of the condition that $H^i(X,\mathcal{O}_X) =0$ for $i\geq 1$, for each ideal sheaf $\mathcal{I}$ with corresponding point $[\mathcal{I}]\in \text{Hilb}^n(X)$, there are isomorphisms
		$$H^0(E_\bullet^\vee\otimes_{\mathcal{O}_{\text{Hilb}^n(X)}}k([\mathcal{I}]))\to \text{Ext}^1(\mathcal{I},\mathcal{I})$$ and 
		$$H^{-1}(E_\bullet^\vee\otimes_{\mathcal{O}_{\text{Hilb}^n(X)}}k([\mathcal{I}]))\to \text{Ext}^2(\mathcal{I},\mathcal{I}).$$ Here $k([\mathcal{I}])$ denotes the residue field at $[\mathcal{I}]$. 
	\end{rem} 
	
	\subsection{Orientation} \label{section orientation}
	Let $X$ be a smooth projective threefold together with an $N_S$-action, let $n\geq 1$, and let $E_\bullet$ be the perfect obstruction theory from the previous section. Let $K_X$ be the canonical sheaf of $X$. 
	
	If there is an isomorphism $K_X\cong L^{\otimes 2}$, one can construct an orientation on $\text{Hilb}^n(X)$. In the case of a Calabi-Yau threefold, this is due to Y. Toda  \cite[Proposition 3.1]{Toda}. Following email correspondence with Toda describing his result and method of proof, M. Levine handled the case of arbitrary $X$ in the paper \cite{LevineOST}. We quote the result here for later use. 
	\begin{prop}[\cite{LevineOST}, Theorem 6.3] \label{Marcs statement isomorphism determinant to square}
		Suppose that there is an $N_S$-linearized very ample line bundle $L$ on $X$ and that we are given an $N_S$-linearized isomorphism $K_X\to L^{\otimes 2}$. Let $\mathcal{I}$ be the ideal sheaf corresponding to the universal subscheme $i : Z \to \text{Hilb}^n(X)\times X$ of $\text{Hilb}^n(X)\times X$. Let $p_Z:Z\to \text{Hilb}^n(X)$ be the natural projection of $Z$ and let $p_2: \text{Hilb}^n(X)\times X\to X$ be the projection to~$X$. Then $E_{\bullet}$ has virtual rank zero and there exists a canonical $N_S$-equivariant isomorphism 
		\begin{align*} 
			\rho: \det(E_\bullet) \to &(\det(p_{Z,*}\mathcal{O}_Z)\otimes N_{Z/\text{Hilb}^n(X)}(i^*p_2^*L))^{\otimes -2}.
		\end{align*} 
		In particular, the determinant of the perfect obstruction theory is a square if $K_X$ is a square. Here, $N_{Z/\text{Hilb}^n(X)}:\text{Pic}(Z)\to \text{Pic}(\text{Hilb}^n(X))$ is the norm  map, see \cite[Section 7]{DeligneDC}. 
	\end{prop} 
	
	\section{Motivic Donaldson-Thomas invariants of $\mathbb{P}^3$}\label{section definition DT invariant}
	For the remainder of the paper, we take the base-field to be $\mathbb{R}$. In this section, we apply the theory set up in the previous section to define motivic Donaldson-Thomas invariants of $\mathbb{P}^3$. We also do this in the $N_S$-equivariant setting. 
	
	\subsection{Action of $N_S$ on $\mathbb{P}^3$} 
	Let $a,b\in\mathbb{Z}$ be odd and such that $$a,b,3a-b, 3b-a, 3a+b, 3b+a, a-b\text{ and }a+b$$ are nonzero (the denominators in the computations in Section \ref{section computations} are the reason for this assumption). Furthermore, assume that $a> 5b$ (so that we can decide for all the terms we will see in Section \ref{section computations} whether they are positive or negative using Proposition \ref{proposition signs}). 
	
	Consider the action of $N_S$ on $\mathbb{P}^3$ given by
	\begin{align*}
		\left(\begin{matrix} 
			t & 0\\
			0 & t^{-1}
		\end{matrix} \right) &\cdot[X_0:X_1:X_2:X_3] = [t^aX_0: t^{-a}X_1: t^bX_2: t^{-b}X_3]\\
		\sigma &\cdot [X_0:X_1:X_2:X_3]= [-X_1:X_0:-X_3:X_2].
	\end{align*}
	\begin{rem} 
		This action does not have fixed points, but there are the two fixed couples $\{[1:0:0:0], [0:1:0:0] \}$ and $\{[0:0:1:0], [0:0:0:1] \}$.
	\end{rem} 
	\begin{rem} 
		We note that for even $a$ or $b$, the above action would not be well defined, as the map $f:N_S\to \text{SL}_2$ given by 
		$$\left(\begin{matrix} 
			t & 0\\
			0 & t^{-1}
		\end{matrix} \right) \mapsto \left(\begin{matrix} 
			t^a & 0\\
			0 & t^{-a}
		\end{matrix}\right)$$ and $\sigma\to\sigma$ is not a morphism. Indeed, if $f$ were a morphism, we would have $f(\sigma^2) = f(\sigma)^2 = \sigma^2 = -\text{Id}$. But $\sigma^2 = -\text{Id}$, so that $f(\sigma^2) = f(-\text{Id}) = \text{Id}$, which is a contradiction. We would therefore have to send $\sigma$ to 
		$$\left(\begin{matrix} 
			0 & 1\\
			1 & 0
		\end{matrix}\right)$$ but this has determinant $-1$ and is therefore not in $\text{SL}_2$. 
	\end{rem} 
	
	\subsection{Motivic Donaldson-Thomas invariants for $\mathbb{P}^3$} 
	We have the perfect obstruction theory from Section \ref{section POT} on $\text{Hilb}^n(\mathbb{P}^3)$ for $n\geq 1$, which gives rise to one in the $N_S$-equivariant case. Using the generator 
	\[
	\Omega:=X_0dX_1dX_2dX_3-X_1dX_0dX_2dX_3+X_2dX_0dX_1dX_3-X_3dX_0dX_1dX_2
	\]
	for $\omega_{\mathbb{P}^3}(4)\cong\mathcal{O}_{\mathbb{P}^3}$, and the canonical isomorphism 
	$ \mathcal{O}(-4)\cong  \mathcal{O}(-2)^{\otimes 2}$, we have a canonical isomorphism 
	$\omega_{\mathbb{P}^3}(4)\cong \mathcal{O}(-2)^{\otimes 2}$, giving us a well defined orientation on $\text{Hilb}^n(\mathbb{P}^3)$, also in the $N_S$-equivariant case, by Theorem \ref{Marcs statement isomorphism determinant to square}. 
	
	We can therefore make the following definition. 
	\begin{defn} 
		For $n\geq 0$, we let $\tilde{I}_n\in \mathcal{W}(\mathbb{R})\cong\mathbb{Z}$ be the degree of the motivic virtual fundamental class of $\text{Hilb}^n(\mathbb{P}^3)$. 
	\end{defn} 
	We will also consider the following equivariant version.
	\begin{defn} 
		Let $\tilde{I}_n^{N_S}\in H^0(BN_S,\mathcal{W})$ be the degree of the  equivariant virtual fundamental class of $\text{Hilb}^n(\mathbb{P}^3)$. 
	\end{defn} 
	\begin{rem} 
		The structure morphism $p: BN_S\to \text{Spec}(\mathbb{R})$ induces a map $p^*: W(\mathbb{R})=H^0(\text{Spec}(\mathbb{R}),\mathcal{W})\to H^0(BN_S,\mathcal{W})$ and we have that $\tilde{I}_n^{N_S} = p^*\tilde{I}_n$.
	\end{rem} 
	
	\section{Applying virtual localization} 
	Let $n\geq 0$, let $\tilde{I}_n$ be the motivic Donaldson-Thomas invariant for $\mathbb{P}^3$ defined in the previous section, and let $\tilde{I}_n^{N_S}$ be the equivariant version. 
	
	Using \cite[Theorem 6.7]{LevineVLEWC}, if the torus in $N_S$ acts with isolated fixed points, we have that
	\begin{align}\label{formula for tildeI in terms of Euler classes}
		\tilde{I}_n^{N_S} &= \sum_{[\mathcal{I}]\in \text{Hilb}^n(\mathbb{P}^3)^{N_S}} \widetilde{\deg}\left(\frac{e(\text{Ext}^2(\mathcal{I},\mathcal{I}))}{e(\text{Ext}^1(\mathcal{I},\mathcal{I}))}\right).
	\end{align}
	Here, the identity takes place in a suitable localization of $H^0(BN_S,\mathcal{W})$, the sum is over all classes of ideal sheaves $[\mathcal{I}]\in \text{Hilb}^n(\mathbb{P}^3)$ which are fixed by the induced $N_S$-action, $e$ denotes the equivariant Euler class in cohomology of Witt sheaves, and $\widetilde{\deg}$ is the equivariant quadratic degree map. In all our computations, each fixed point will be a copy of $\text{Spec}(\mathbb{R})$, so the degree map is just an identity map and can be ignored. Because any fixed ideal sheaf needs to be invariant under the $\sigma$-action, all fixed ideal sheaves of the $N_S$-action are of even degree. From this and (\ref{formula for tildeI in terms of Euler classes}) we immediately find the following. 
	\begin{cor} 
		$\tilde{I}_n = 0 $ whenever $n$ is odd.
	\end{cor} 
	If $n$ is even, we need to compute the trace of the virtual representation 
	\begin{equation}\label{representation} 
		\text{Ext}^2(\mathcal{I},\mathcal{I}) - \text{Ext}^1(\mathcal{I},\mathcal{I})
	\end{equation} 
	of $N_S$ for all classes of ideal sheaves $\mathcal{I}$ that are isolated fixed points of the $N_S$-action; if the $N_S$-fixed locus on $\text{Hilb}^n(X)$ has non-isolated fixed points, the method requires more work, and will not be discussed further here. We can compute the equivariant Euler classes from this trace using Theorem \ref{Marcs computation of euler classes}, which we will apply in Section \ref{section computations} to compute $\tilde{I}_n$ for $n\leq 6$. 
	
	In this section, we first study how the orientation from Proposition \ref{Marcs statement isomorphism determinant to square} induces signs on the resulting degrees of virtual fundamental classes, see Section \ref{section signs}. We then compute the trace of the virtual representation (\ref{representation}) using the methods of \cite{MaulikI} and \cite{MaulikII}, see Section \ref{section trace computation}. 
	\begin{nota} 
		Throughout this section, we denote $U_i = \{X_i\neq 0\}\subset\mathbb{P}^3$ for $i\in\{0,\cdots, 3\}$.  
	\end{nota}
	
	\subsection{Signs coming from the orientation}\label{section signs}
	We now study the signs that the orientation from Proposition \ref{Marcs statement isomorphism determinant to square} induces on the equivariant virtual fundamental classes $\tilde{I}_n^{N_S}$. 
	\begin{prop}\label{proposition signs}
		The orientation on the perfect obstruction theory on $\text{Hilb}^n(\mathbb{P}^3)$ from Proposition \ref{Marcs statement isomorphism determinant to square} gives rise to an oriented basis of $\text{Ext}^2(\mathcal{I},\mathcal{I}) - \text{Ext}^1(\mathcal{I},\mathcal{I})$ for each $N_S$-fixed ideal sheaf $\mathcal{I}$. Every even negative weight in the trace induces a minus sign to the corresponding Euler class.   
	\end{prop}
	\begin{proof}
		Let $\mathcal{I}\subset\mathcal{O}_{\mathbb{P}^3}$ be an ideal sheaf corresponding to an $N_S$-fixed point. Let $I$ be the image of $\mathcal{I}$ on $U_0$. 
		
		Noting that $n$ is even, $n=2m$, we have the oriented perfect obstruction theory on $\text{Hilb}^m(\mathbb{P}^3)$. Because this perfect obstruction theory has virtual rank zero, we have that $\text{Ext}^1(I,I)$ and $\text{Ext}^2(I,I)$ have the same rank $r$. We let $e_1,\cdots, e_r$,  $f_1,\cdots, f_r$ be bases for $\text{Ext}^2(I,I)$, resp. $\text{Ext}^1(I,I)$, that yield an oriented basis for the virtual vector space $\text{Ext}^2(I,I) - \text{Ext}^1(I,I)$, that is $$\rho_m(\det(\text{Ext}^2(I,I) - \text{Ext}^1(I,I))) := \rho_m((e_1\wedge\cdots\wedge e_r) \otimes (f_1\wedge\cdots\wedge f_r)^{-1}) = 1.$$
		
		By \cite[Theorem 6.4]{LevineOST}, $\rho_m$ is $N_S$-equivariant and so we have that $$(-\sigma)(e_1),\cdots,(-\sigma)(e_r); (-\sigma)(f_1),\cdots, (-\sigma)(f_r)$$ yields an oriented basis at $\sigma(I)$. Looking at the construction of the orientation for the perfect obstruction theory on $\text{Hilb}^n(\mathbb{P}^3)$, as given by \cite[Theorem 6.3]{LevineOST},  it follows that 
		\[
		e_1,\cdots, e_r, (-\sigma)(e_1),\cdots,(-\sigma)(e_r); f_1,\cdots, f_r, (-\sigma)(f_1),\cdots, (-\sigma)(f_r)	
		\]
		is a pair of bases for 	 $\text{Ext}^2(\mathcal{I},\mathcal{I})$, resp. $\text{Ext}^1(\mathcal{I},\mathcal{I})$ that together yield an oriented basis for $\text{Ext}^2(\mathcal{I},\mathcal{I}) - \text{Ext}^1(\mathcal{I},\mathcal{I})$ with respect to the orientation $\rho_n$. If we now change each of these bases by the same permutation of $\{1,\ldots, 2r\}$, each of the determinants will change by the same sign, so the resulting bases will still yield an oriented basis for $\text{Ext}^2(\mathcal{I},\mathcal{I}) - \text{Ext}^1(\mathcal{I},\mathcal{I})$.
		
		Thus the pair of bases
		$$e_1, (-\sigma)(e_1), \cdots, e_r, (-\sigma)(e_r); f_1, (-\sigma)(f_1), \cdots, f_r, (-\sigma)(f_r)$$ for $\text{Ext}^2(\mathcal{I},\mathcal{I})$, resp. $\text{Ext}^1(\mathcal{I},\mathcal{I})$ gives a basis 
		for $\text{Ext}^2(\mathcal{I},\mathcal{I}) - \text{Ext}^1(\mathcal{I},\mathcal{I})$ on $U_0$ and $U_1$, which is compatible with the relative orientation $\rho_n$ as given by Proposition \ref{Marcs statement isomorphism determinant to square}. In order to use Proposition \ref{Marcs computation of euler classes} to compute Euler classes, we need to follow Remark \ref{remark sign Euler classes} and use a basis which is compatible with the choice described there. This means: for each $e_i$ with negative weight, one has to switch $e_i$ and $(-\sigma)(e_i)$ in the above basis, because the positive weight always has to be the first; we then  replace $(-\sigma)(e_i)$ with $\sigma(e_i)$ . The same holds for the $f_i$. If the weight of a negative $e_i$ or $f_i$ is odd, these switches do not contribute a sign change to the Euler class, as $(-\sigma)(e_i) = -\sigma(e_i)$ in this case, so 
		\[
		e_i\wedge(-\sigma)(e_i)= e_i\wedge(-\sigma(e_i))=\sigma(e_i)\wedge e_i, 
		\]
		and similarly for $f_i$.
		If the weight of a negative $e_i$ or $f_i$ is even,  this contributes a sign to the Euler class since then $(-\sigma)(e_i) = \sigma(e_i)$, $(-\sigma)(f_i) = \sigma(f_i)$. 
		
		One can repeat this construction on $U_2$ and $U_3$. This proves the desired statement. 
	\end{proof} 
	
	\subsection{Computing the trace}\label{section trace computation}
	In order to find the trace of the virtual representation (\ref{representation}), we use the strategy of \cite{MaulikI} and \cite{MaulikII}. The proofs are almost the same as in \cite[Section 4]{MaulikI}, but included here for the reader's convenience. We first show the following. 
	\begin{lem}[See Formula (9) of \cite{MaulikI}]\label{lemma Extgroups and Cech cohomology}
		For an ideal sheaf $[\mathcal{I}]\in \text{Hilb}^n(\mathbb{P}^3)$, fixed under the $N_S$-action, we have that
		\begin{align*} 
			\text{Ext}^1(\mathcal{I},\mathcal{I}) - \text{Ext}^2(\mathcal{I},\mathcal{I}) &= \sum_{i=0}^3 \left(H^0(U_i, \mathcal{O}_{\mathbb{P}^3}) - \sum_{j=0}^3 (-1)^jH^0(U_i,\text{Ext}^j(\mathcal{I},\mathcal{I}))\right)
		\end{align*} 
		as virtual representations of $N_S$. 
	\end{lem} 
	\begin{rem}
		Note that one can view the above as an $N_S$-representation by considering it as the sum of 
		\begin{align*} 
			&\left(H^0(U_0, \mathcal{O}_{\mathbb{P}^3}) - \sum_{j=0}^3 (-1)^jH^0(U_0,\text{Ext}^j(\mathcal{I},\mathcal{I}))\right)\\
			\quad &\oplus \left(H^0(U_1, \mathcal{O}_{\mathbb{P}^3}) - \sum_{j=0}^3 (-1)^jH^0(U_1,\text{Ext}^j(\mathcal{I},\mathcal{I}))\right)\end{align*}  and the corresponding term for $U_2$ and $U_3$. 
	\end{rem}
	To prove the lemma, we start by finding a helpful expression for $\text{Ext}^2(\mathcal{I},\mathcal{I}) - \text{Ext}^1(\mathcal{I},\mathcal{I})$ in terms of equivariant Euler characteristics. 
	\begin{constr}
		For two coherent sheaves $\mathcal{F}$ and $\mathcal{G}$ on $\mathbb{P}^3\times\text{Hilb}^n(\mathbb{P}^3)$ with $\mathcal{F}$ flat over $\text{Hilb}^n(\mathbb{P}^3)$, we have the derived Hom-set $R\text{Hom}(\mathcal{F},\mathcal{G})$ in the bounded derived category $D^b(\mathbb{R})$. Its Euler characteristic is given by $$\tilde{\chi}(\mathcal{F},\mathcal{G}) =  \sum_{i=0}^3(-1)^i\dim_{\mathbb{R}}(\text{Ext}^i(\mathcal{F},\mathcal{G}))\in \mathbb{Z}.$$ 
	\end{constr}
	For two coherent sheaves $\mathcal{F}$ and $\mathcal{G}$ on $\mathbb{P}^3\times\text{Hilb}^n(\mathbb{P}^3)$, an $N_S$-linearization on $\mathcal{F}$ and $\mathcal{G}$ defines an $N_S$-action on $\text{Ext}^i(\mathcal{F},\mathcal{G})$, which is therefore an $N_S$-representation. This yields a refined Euler characteristic. 
	\begin{defn}
		Let $\mathcal{F}$ and $\mathcal{G}$ be coherent sheaves on $\mathbb{P}^3$, each endowed with an $N_S$-linearization. The \textit{Euler characteristic} $\chi(\mathcal{F},\mathcal{G})$ is given by the alternating sum of virtual $N_S$-representations $$\chi(\mathcal{F},\mathcal{G}) = \sum_{i=0}^3(-1)^i \text{Ext}^i(\mathcal{F},\mathcal{G})\in K_0(N_S-\text{Reps}).$$
	\end{defn}	
	We now make the following observation. 
	\begin{lem} \label{exts in terms of chi}
		Let $\mathcal{I}$ be an ideal sheaf on $\mathbb{P}^3$ which is stable under the $N_S$-action. We have that 
		$$\chi(\mathcal{I},\mathcal{I}) - \chi(\mathcal{O}_{\mathbb{P}^3}, \mathcal{O}_{\mathbb{P}^3}) = \text{Ext}^2(\mathcal{I},\mathcal{I}) - \text{Ext}^1(\mathcal{I},\mathcal{I})$$
		in $K_0(N_S-\text{Reps})$.	\end{lem} 
	\begin{proof}
		Writing out gives 
		\begin{align*} 
			\chi(\mathcal{I},\mathcal{I}) - \chi(\mathcal{O}_{\mathbb{P}^3}, \mathcal{O}_{\mathbb{P}^3}) &= \sum_{i=0}^3(-1)^i\text{Ext}^i(\mathcal{I},\mathcal{I}) - \text{Hom}(\mathcal{O}_{\mathbb{P}^3},\mathcal{O}_{\mathbb{P}^3})\\
			&= -\text{Ext}^3(\mathcal{I},\mathcal{I}) + \text{Ext}^2(\mathcal{I},\mathcal{I}) - \text{Ext}^1(\mathcal{I},\mathcal{I})\\
			&\quad  + \text{Hom}(\mathcal{I},\mathcal{I}) - \text{Hom}(\mathcal{O}_{\mathbb{P}^3},\mathcal{O}_{\mathbb{P}^3})\\
			&= \text{Ext}^2(\mathcal{I},\mathcal{I}) - \text{Ext}^1(\mathcal{I},\mathcal{I})
		\end{align*} 
		because $\text{Ext}^3(\mathcal{I},\mathcal{I}) = 0$ by \cite[Lemma 2]{MaulikI} and $\text{Hom}(\mathcal{I},\mathcal{I}) = \mathcal{O}_{\mathbb{P}^3}$. 
	\end{proof}
	\begin{proof}[Proof of Lemma \ref{lemma Extgroups and Cech cohomology}]
		Consider the local-to-global spectral sequence 
		$$E_2^{p,q} = H^p(\mathbb{P}^3, \mathcal{E}xt^q(\mathcal{I},\mathcal{I})) \implies \text{Ext}^{p+q}(\mathcal{I},\mathcal{I}).$$
		Note that $E_2^{p,q} = 0$ if $p > 3$ or $q>3$. Note that 
		\begin{align*} 
			\sum_{p,q}(-1)^{p+q}E_s^{p,q} &= \sum_{p,q}(-1)^{p+q}(\ker(d_s^{p,q}) + \text{im}(d_s^{p,q}))\\
			&= \sum_{p,q}(-1)^{p+q}E_{s+1}^{p,q}
		\end{align*} 
		for all $s\geq 2$. In particular, this sums the infinity pages with the right signs, so that we find
		\begin{align*}
			\chi(\mathcal{I},\mathcal{I}) &= \sum_{i,j =0}^3(-1)^{i+j}H^i(\mathbb{P}^3, \mathcal{E}xt^j(\mathcal{I},\mathcal{I})).
		\end{align*}
		Using the $U_i$ as a cover of $\mathbb{P}^3$, we can compute the above cohomology groups by computing \v{C}ech cohomology. Note that as $\mathcal{I}$ is only supported on points $[x_0,\ldots, x_3]\in \mathbb{P}^3$ with exactly one $x_i\neq0$, we have that $\mathcal{I} = \mathcal{O}_{\mathbb{P}^3}$ on the intersection of two or more $U_i$. Therefore, the \v{C}ech complex for $\mathcal{E}xt^j(\mathcal{I},\mathcal{I})$ is concentrated in degree zero, showing that
		\[ 
		\text{Ext}^j(\mathcal{I},\mathcal{I})=\oplus_iH^0(U_i, \mathcal{E}xt^j(\mathcal{I},\mathcal{I})).
		\]
		By Lemma \ref{exts in terms of chi}, this gives the desired statement. 
	\end{proof} 
	
	For a representation $\rho$ of $\mathbb{G}_m$ on an $\mathbb{R}$-vector space $V$, we let $\text{tr}_V(t)$ be the trace of $\rho$, i.e., if $V=\oplus_mV_m$ where $\rho(t)(v)=t^mv$ for $v\in V_m$, then $\text{tr}_V(t):=\sum_m \dim_{\mathbb{R}}V_m\cdot t^m$; we will only consider  $\text{tr}_V(t)$ for $V$ such that $\dim_{\mathbb{R}}V_m$ is finite for each $m$.
	
	\begin{prop}[See Formula 12 of \cite{MaulikI}]\label{formula for trace of our representation}
		Let $\mathcal{I}$ be an $N_S$-stable ideal sheaf. For $i\in\{0,\cdots, 3\}$, write $R = \mathbb{R}[x,y,z] \cong H^0(U_i, \mathcal{O}_{\mathbb{P}^3})$, with $x,y,z$ the standard coordinates $X_j/X_i$, $j\neq i$, and let $I$ be the image of the ideal sheaf $\mathcal{I}$. Let $\pi_I = \{(i,j,k): x^iy^jz^k\notin I \}$ and suppose that $s_1,s_2,s_3$ are the tangent weights on $U_i$, i.e. $t\cdot (x,y,z) = (s_1x,s_2y,s_3z)$. Set $$Q_i(t) = \text{tr}_{R/I}(t) = \sum_{(i,j,k)\in \pi_I}s_1^is_2^js_3^k.$$ 
		Then 
		\begin{align*} 
			\text{tr}_{\text{Ext}^1(I,I)-\text{Ext}^2(I,I)}(t) &= \frac{s_1s_2s_3Q(t) - Q(t^{-1}) + Q(t)Q(t^{-1})(1-s_1)(1-s_2)(1-s_3)}{s_1s_2s_3}.
		\end{align*}  
	\end{prop}
	Adding the traces on different $U_i$ and filling in the correct tangent weights gives the trace of the virtual tangent space of $\mathcal{I}$. 
	\begin{proof}[Proof of Proposition \ref{formula for trace of our representation}]
		We have that 
		\begin{align*}
			\text{tr}_R(t) &= \sum_{i,j,k} s_1^is_2^js_3^k\\
			&= \left(\sum_i s_1^i\right)\left(\sum_j s_2^j\right)\left(\sum_k s_3^k\right)\\
			&=  \frac{1}{(1-s_1)(1-s_2)(1-s_3)}
		\end{align*}
		Consider a resolution of the ideal $I$ given by 
		\begin{equation} \label{resolution of the ideal}
			0\to F_r\to\cdots\to F_1\to I\to 0
		\end{equation} 
		such that each term is of the form $F_i = \bigoplus_j R(d_{ij})$ for $d_{ij} = (d_{ij}^1, d_{ij}^2, d_{ij}^3)\in\mathbb{Z}^3$. Here $R(d^1, d^2, d^3)$ is the representation with 
		\[
		t\cdot x^iy^jz^k=s_1^{i+d^1}s_2^{j+d^2}s_3^{k+d^3}x^iy^jz^k.
		\]
		One can always form a resolution like this if the ideal $I$ is homogeneous: if $I = (f_0,\cdots, f_m)$ where $f_i$ has degree $d_i$ then we can take $F_1 = \bigoplus_{i=0}^mR(d_i)$, define $F_2$ based on the relations between the $f_i$ and continue like that until we have a resolution. Consider the corresponding Poincar\'e polynomial 
		$$P(t) = \sum_{i,j}(-1)^is_1^{d_{ij}^1}s_2^{d_{ij}^2}s_3^{d_{ij}^3}.$$ 
		Note that for any $d = (d^1,d^2,d^3)\in\mathbb{Z}^3$, we have that $$\text{tr}_{R(d)}(t) = \sum_{i,j,k} t\cdot x^{i+d^1}y^{j+d^2}z^{k+d^3} = s_1^{d^1}s_2^{d^2}s_3^{d^3}\text{tr}_R(t).$$ This implies that 
		\begin{align*}
			\text{tr}_{F_*}(t) &= \sum_{i,j}(-1)^i\text{tr}_{R(d_{ij})}(t)\\
			&= \text{tr}_R(t)\sum_{i,j}(-1)^is_1^{d_{ij}^1}s_2^{d_{ij}^2}s_3^{d_{ij}^3} \\
			&= \frac{P(t)}{(1-s_1)(1-s_2)(1-s_3)}
		\end{align*}
		We have the exact sequence 
		$$0\to F_*\to R\to R/I\to 0$$ and so 
		$\text{tr}_R(t) = \text{tr}_{F_*}(t) + \text{tr}_{R/I}(t)$.
		This implies that
		$$Q(t) = \frac{1 - P(t)}{(1-s_1)(1-s_2)(1-s_3)}.$$
		We use the resolution $F_*$ to compute that 
		\begin{align*}
			\chi(I,I) &= \sum_{i,j,k,l}(-1)^{i+k}\text{Hom}(R(d_{ij}), R(d_{kl}))\\
			&= \sum_{i,j,k,l}(-1)^{i+k}R(d_{kl} - d_{ij})
		\end{align*}
		This tells us that the trace is 
		\begin{align*}
			\text{tr}_{\chi(I,I)}(t) &= \sum_{i,j,k,l}(-1)^{i+k}s_1^{d_{kl}^1-d_{ij}^1}s_2^{d_{kl}^2-d_{ij}^2}s_3^{d_{kl}^3-d_{ij}^3}\text{tr}_R(t)\\
			&= \frac{P(t^{-1})P(t)}{(1-s_1)(1-s_2)(1-s_3)}
		\end{align*}
		and so 
		\begin{align*}
			\text{tr}_{R-\chi(I,I)}(t) &= \frac{1-P(t^{-1})P(t)}{(1-s_1)(1-s_2)(1-s_3)}.
		\end{align*}
		We have that $P(t) = -(1-s_1)(1-s_2)(1-s_3)Q(t) + 1$ and 
		\begin{align*} 
			P(t^{-1}) &= -(1-s_1^{-1})(1-s_2^{-1})(1-s_3^{-1})Q(t^{-1}) + 1\\
			&= s_1s_2s_3(1-s_1)(1-s_2)(1-s_3)Q(t^{-1}) + 1
		\end{align*} 
		and so we can rewrite this as
		\begin{align*}
			\text{tr}_{R - \chi(I,I)}(t) &= \frac{s_1s_2s_3Q(t) - Q(t^{-1}) + Q(t)Q(t^{-1})(1-s_1)(1-s_2)(1-s_3)}{s_1s_2s_3}.
		\end{align*}
		as desired. 
	\end{proof} 
	
	\section{Computations of $\tilde{I}_n$ for  $n\leq 6$}\label{section computations}
	For $n\geq 0$, let $\tilde{I}_n$ be the motivic Donaldson-Thomas invariant for $\mathbb{P}^3$ defined in Section \ref{section definition DT invariant}. For $i\in\{0,\cdots, 3\}$, we let $U_i = \{X_i\neq 0\}\subset\mathbb{P}^3$. 
	
	We use the virtual localization formula (\ref{formula for tildeI in terms of Euler classes}) together with Proposition \ref{Marcs computation of euler classes}, Lemma \ref{lemma Extgroups and Cech cohomology} and Proposition \ref{formula for trace of our representation} to compute the motivic Donaldson-Thomas invariants $\tilde{I}_n\in H^*(\text{Spec}(\mathbb{R}),\mathcal{W})\cong\mathbb{Z}$ for $n\leq 6$, which will prove Theorem \ref{main theorem}.  
	\begin{constr}
		Throughout, we use coordinates 
		\begin{itemize} 
			\item $x = \frac{X_1}{X_0}, y = \frac{X_2}{X_0}, z = \frac{X_3}{X_0}$ on $U_0$.
			\item $u = \frac{X_0}{X_1}, v = \frac{X_2}{X_1}, w = \frac{X_3}{X_1}$ on $U_1$.
			\item $x' = \frac{X_0}{X_2}, y' = \frac{X_1}{X_2}, z' = \frac{X_3}{X_2}$ on $U_2$.
			\item $u' = \frac{X_0}{X_3}, v' = \frac{X_1}{X_3}, w' = \frac{X_2}{X_3}$ on $U_3$. 
		\end{itemize} 
	\end{constr}
	\begin{constr} 
		As $[t^aX_0: t^{-a}X_1: t^bX_2: t^{-b}X_3] = [X_0: t^{-2a}X_1: t^{b-a}X_2: t^{-a-b}X_3]$, we have that on $U_0$
		$$\left(\begin{matrix} 
			t & 0 \\ 0 & t^{-1} \end{matrix}\right) \cdot (x,y,z) = (s_1x, s_2y, s_3z)$$ where $s_1 = t^{-2a}, s_2 = t^{b-a}, s_3 = t^{-a-b}$. Repeating this for $U_1,U_2$ and $U_3$, we see that the ``tangent weights" of the action are 
		\begin{align*}
			U_0: &\quad s_1 = t^{-2a}, s_2 = t^{b-a}, s_3 = t^{-a-b} \\
			U_1: &\quad s_1 = t^{2a}, s_2 = t^{b+a}, s_3 = t^{a-b}\\
			U_2: &\quad s_1 = t^{a-b}, s_2 = t^{-b-a}, s_3 = t^{-2b}\\
			U_3: &\quad s_1 = t^{a+b}, s_2 = t^{b-a}, s_3 = t^{2b}
		\end{align*}
	\end{constr} 
	\begin{rem} 
		Note the symmetries between those tangent weights, which imply that for an $N_S$-stable ideal sheaf $\mathcal{I}$ of even degree $n$ supported on $\{[1:0:0:0], [0:1:0:0] \}$, we can compute the trace on $U_1$ as $V(t^{-1})$, where $V(t)$ is the trace on $U_0$. From $V(t) + V(t^{-1})$, we can find the Euler class $e(a,b)$ corresponding to $\mathcal{I}$ by Lemma \ref{lemma Extgroups and Cech cohomology}. Furthermore, for each $N_S$-stable ideal sheaf $\mathcal{I}$ suppported on $\{[1:0:0:0], [0:1:0:0] \}$, we have a corresponding $N_S$-stable ideal sheaf $\mathcal{I}'$ supported on $\{[0:0:1:0], [0:0:0:1] \}$. Using the symmetries between the tangent weights again, we see that $\mathcal{I}'$ has Euler class $e(b,a)$. 
	\end{rem}
	\begin{rem}
		By Proposition \ref{proposition signs}, each even negative weight in the trace induces a minus sign to the corresponding Euler class. For all Euler classes computed in this section, the signs introduced by these even negative weights have been taken into account in the computation, and will not be mentioned further. 
	\end{rem}
	
	In the following discussion, we will say that an ideal sheaf $\mathcal{I}\subset\mathcal{O}_{\mathbb{P}^3}$ is supported on some subset $S\subset \mathbb{P}^3$ if the corresponding closed scheme $\text{Spec}(\mathcal{O}_{\mathbb{P}^3}/\mathcal{I})$ has support contained in $S$. We will often identify an $N_S$-stable ideal sheaf with the corresponding $N_S$-fixed point in $\text{Hilb}(\mathbb{P}^3)$, referring to both as a fixed point.
	
	\subsection{The computation for $n=2$}
	Note that there are two $N_S$-stable ideal sheaves of length two: 
	\begin{enumerate} 
		\item The subscheme supported on $\{[1:0:0:0], [0:1:0:0] \}$ of which the ideal sheaf $\mathcal{I}$ is given by $(x,y,z)\subset \mathbb{R}[x,y,z]$ on $U_0$ and by $(u,v,w)\subset \mathbb{R}[u,v,w]$ on $U_1$. 
		\item The subscheme supported on $\{[0:0:1:0], [0:0:0:1] \}$, again with the ideal given by $(x',y',z')\subset \mathbb{R}[x',y',z']$ on $U_2$ and by $(u',v',w')\subset \mathbb{R}[u',v',w']$ on $U_3$. 
	\end{enumerate} 
	\begin{prop} \label{computation n=2}
		The Euler class corresponding to the point (1) is 
		$$e_{11} = \frac{(3a - b)(3a+b)}{(a-b)(a+b)}$$ and the class corresponding to (2) is 
		$$e_{12} = \frac{(3b - a)(3b+a)}{(b-a)(a+b)}.$$
		We have that $\tilde{I}_2 = e_{11} + e_{12} = 10$. 
	\end{prop}  
	\begin{proof} 
		We start with the point (1). Plugging $Q= 1$ into the formula from Proposition \ref{formula for trace of our representation} yields that the trace of $\text{Ext}^1(\mathcal{I},\mathcal{I}) - \text{Ext}^2(\mathcal{I},\mathcal{I})$ is given by
		$$V_2 = s_1^{-1} + s_2^{-1} + s_3^{-1} - s_1^{-1}s_2^{-1} - s_1^{-1}s_3^{-1} - s_2^{-1}s_3^{-1}.$$
		Filling in the tangent weights for $U_0$ and $U_1$ and adding up gives 
		\begin{align*} 
			(t^{a-b} + t^{b-a}) + (t^{a+b} + t^{-a-b}) &- ((t^{3a-b} + t^{-3a+b}) + (t^{3a+b} + t^{-3a-b}).
		\end{align*}
		Now using Proposition \ref{Marcs computation of euler classes}, we see that this leads to an Euler class 
		$$e_{11} = \frac{(3a - b)(3a+b)}{(a-b)(a+b)}.$$
		We note that for $a,b$ both congruent to $1$ modulo $4$, we have, modulo $4$, that $3b-a = 2, 3b+a = 0, b-a = 0, a+b = 2$ so that the use of Proposition \ref{Marcs computation of euler classes} gives a plus. One can check that the sign does not change with all other possible congruences of $a,b$ modulo $4$. 
		
		Repeating the process for the second fixed point boils down to switching $a$ and $b$, and gives us the Euler class 
		$$e_{12} = \frac{(3b - a)(3b+a)}{(b-a)(a+b)}. $$
		Adding these we get $e_{11} + e_{12} = 10$. 
	\end{proof} 
	
	\subsection{The computation for $n=4$} 
	Note that the $N_S$-action exchanges $U_0$ and $U_1$, with $N_S$ acting by $x\mapsto -u, y\mapsto w, z\mapsto -v$. From this we see that to be stable under the $N_S$-action, an ideal sheaf supported on $\{[1:0:0:0], [0:1:0:0] \}$ of length four must correspond to one of the following:
	\begin{enumerate} 
		\item $(x^2,y,z)$ on $U_0$ and $(u^2, v,w)$ on $U_1$.
		\item $(x,y^2,z)$ on $U_0$ and $(u, v, w^2)$ on $U_1$.
		\item $(x,y,z^2)$ on $U_0$ and $(u, v^2 ,w)$ on $U_1$.
	\end{enumerate} 
	One can find similar classes for stable ideal sheaves of length four supported on $\{[0:0:1:0], [0:0:0:1] \}$. For $n=4$ there are thus seven fixed points if we add the ideal sheaf supported on all four points. We will prove the following result. 
	\begin{prop}\label{computation n=4}
		The fixed point supported on all four points has Euler class 
		$$e_{21} = -\frac{(3a-b)(3a+b)(3b-a)(3b+a)}{(a-b)^2(a+b)^2}.$$
		The points of type (1) have Euler class zero, and the points in type (2) give rise to Euler classes $$e_{22} = \frac{(3a-b)^2(3a+b)(2a-b)}{b(a-b)^2(a+b)}\text{ and } e_{23} = \frac{(3b-a)^2(3b+a)(2b-a)}{a(b-a)^2(a+b)}.$$
		The points of type (3) have Euler classes 
		$$e_{24} = -\frac{(3a-b)(3a+b)^2(2a+b)}{b(a-b)(a+b)^2}\text{ and } e_{25} = -\frac{(3b-a)(3b+a)^2(2b+a)}{a(b-a)(a+b)^2}.$$
		We have that $\tilde{I}_4 = e_{21} + e_{22} + e_{23} + e_{24} + e_{25} = 25$. 
	\end{prop}
	\begin{proof} 
		First, there is the subscheme which consists of all four points. The Euler class is 
		$$e_{21} = e_{11}e_{12} = -\frac{(3a-b)(3a+b)(3b-a)(3b+a)}{(a-b)^2(a+b)^2}.$$	
		For the fixed point of type (1) supported on $\{[1:0:0:0], [0:1:0:0] \}$, plugging $Q = 1+s_1$ into the formula from Proposition \ref{formula for trace of our representation} gives 
		\begin{align} \label{V4} 
			V_4 &= s_1^{-1} + s_2^{-1} + s_3^{-1} - s_1^{-1}s_2^{-1} - s_1^{-1}s_3^{-1}  - s_2^{-1}s_3^{-1}\\
			&\quad + s_1^{-2} + s_1s_3^{-1} + s_1s_2^{-1} -s_1s_2^{-1}s_3^{-1} - s_1^{-2}s_2^{-1} - s_1^{-2}s_3^{-1}.\nonumber 
		\end{align} 
		Filling in the tangent weights on $U_0$ gives 
		\begin{align*} 
			t^{2a} &+ t^{a-b} + t^{a+b} - t^{3a-b} - t^{2a} - t^{3a+b} + \\
			t^{4a} &+ t^{b-a} + t^{-b-a} - 1 - t^{5a-b} - t^{5a+b}.
		\end{align*} 
		and filling in the tangent weights of $U_1$ gives the same but with $t$ and $t^{-1}$ swapped. This yields an Euler class which is zero, because of the term $1$ which represents a trivial subrepresentation of $\text{Ext}^2(\mathcal{I},\mathcal{I})$, which contributes a zero to the numerator. Similarly, the fixed point supported on $\{[0:0:1:0], [0:0:0:1] \}$ of type~(1) has Euler class equal to zero. 
		
		Now consider the fixed point supported on $\{[1:0:0:0], [0:1:0:0] \}$ of type~(2). Filling in $V_4$ with the tangent weights of $U_0$ where we switch $s_1$ with $s_2$ gives 
		\begin{align*} 
			t^{a-b} &+ t^{2a} + t^{a+b} - t^{3a-b} - t^{2a} - t^{3a+b} + \\
			t^{2a-2b} &+ t^{2b} + t^{b+a} - t^{2a+2b} - t^{3a-b} - t^{4a-2b}.
		\end{align*} 
		If we fill in $V_4$ with the tangent weights of $U_1$ where we switch $s_1$ with $s_3$, this gives the same result with $t$ and $t^{-1}$ switched. This yields the Euler class~$e_{22}$. Similarly, the fixed point supported on $\{[0:0:1:0], [0:0:0:1] \}$ of type (2) has Euler class $e_{23}$. 
		
		Finally, for the fixed point supported on $\{[1:0:0:0], [0:1:0:0] \}$ of type (3), filling in $V_4$ with the tangent weights of $U_0$ where we switch $s_1$ with $s_3$ gives 
		\begin{align*}
			t^{a+b} &+ t^{a-b} + t^{2a} - t^{3a+b} - t^{2a} - t^{3a-b} + \\
			t^{2a+2b} &+ t^{a-b} + t^{-2b} - t^{2a-2b} - t^{4a+2b} - t^{3a+b}.
		\end{align*}
		Again, filling in $V_4$ with the tangent weights of $U_1$ and $s_1$ and $s_2$ switched boils down to switching $t$ and $t^{-1}$. This yields the Euler class $e_{24}$, so that the fixed point supported on $\{[0:0:1:0], [0:0:0:1] \}$ of type (3) has  Euler class $e_{25}$ as desired.
		
		One can check that for all possible congruence classes of $a,b$ modulo $4$, no signs are introduced because of the use of Proposition \ref{Marcs computation of euler classes}. 
	\end{proof} 
	
	\subsection{The computation for $n=6$}
	We will now prove the following. 
	\begin{prop}\label{computation Itilde6} 
		We have that $\tilde{I}_6 = -50$. 
	\end{prop}  
	There are three types of stable ideal sheaves of length six. First, from the lower degrees, we have the Euler classes 
	\begin{align*}
		e_{31} &= e_{11}\cdot e_{23}\\
		e_{32} &= e_{11}\cdot e_{25}\\
		e_{33} &= e_{12}\cdot e_{22}\\
		e_{34} &= e_{12} \cdot e_{24}
	\end{align*}
	Another class of ideal sheaves are the complete intersections, i.e. those supported on $[1:0:0:0]$ and $[0:1:0:0]$ which correspond to either:
	\begin{enumerate}
		\item $(x^3,y,z)$ on $U_0$ and $(u^3,v,w)$ on $U_1$.
		\item $(x,y^3,z)$ on $U_0$ and $(u,v,w^3)$ on $U_1$.
		\item $(x,y,z^3)$ on $U_0$ and $(u,v^3,w)$ on $U_1$.
	\end{enumerate}
	Again, each of those have a corresponding fixed point supported on $[0:0:1:0]$ and $[0:0:0:1]$. 
	\begin{lem} 
		The fixed points of type (1) have Euler class equal to zero. The points of type (2) yield the Euler classes 
		$$e_{35} = \frac{(3a-b)^2(3a+b)(a+3b)(5a-3b)(2a-b)^2}{3b^2(a-b)^3(a+b)(3b-a)}$$
		and
		$$e_{36} = \frac{(3b-a)^2(3b+a)(b+3a)(5b-3a)(2b-a)^2}{3a^2(b-a)^3(a+b)(3a-b)}.$$
		The points of type (3) have Euler classes $$e_{37} =  - \frac{(3a-b)(3a+b)^2(a-3b)(5a+3b)(2a+b)^2}{3b^2(a-b)(a+b)^3(a+3b)}$$ and 
		$$e_{38} =  -\frac{(3b-a)(3b+a)^2(b-3a)(5b+3a)(2b+a)^2}{3a^2(b-a)(a+b)^3(b+3a)}.$$
	\end{lem} 
	\begin{proof} 
		Consider the ideal of type (1). Plugging $Q = 1 + s_1 + s_1^2$ into the formula from Proposition \ref{formula for trace of our representation} gives the trace 
		\begin{align*} 
			V_{61} &= s_1^{-1} + s_2^{-1} + s_3^{-1} - s_1^{-1}s_2^{-1} - s_1^{-1}s_3^{-1} - s_2^{-1}s_3^{-1} \\ 
			&\quad + s_1^{-2} + s_1s_2^{-1} + s_1s_3^{-1} - s_1s_2^{-1}s_3^{-1} - s_1^{-2}s_2^{-1} - s_1^{-2}s_3^{-1}   \\
			&\quad + s_1^{-3} + s_1^2s_2^{-1} + s_1^2s_3^{-1} - s_1^2s_2^{-1}s_3^{-1} - s_1^{-3}s_2^{-1} - s_1^{-3}s_3^{-1}.
		\end{align*}
		Note that the first two lines are precisely the trace (\ref{V4}) we had in the $n=4$ case. Filling in the tangent weights of $U_0$ in the last line of $V_{61}$ gives 
		\begin{align*}
			t^{6a} &+ t^{-3a + b} + t^{-3a - b} - t^{-2a} - t^{7a + b} - t^{7a-b}.
		\end{align*}
		There is nothing here which cancels the trivial subrepresentation coming from the second line, implying that the corresponding Euler class is zero. Similarly, the Euler class of the corresponding fixed point supported on $[0:0:1:0]$ and $[0:0:0:1]$ is zero. 
		
		For the fixed point supported on $[1:0:0:0]$ and $[0:1:0:0]$ of type (2), filling in the tangent weights of $U_0$ in the last line of $V_{61}$ and switching $s_1$ and $s_2$ gives 
		\begin{align*}
			t^{3a-3b} &+ t^{-a + 3b} + t^{2b} - t^{3b+a} - t^{5a-3b} - t^{4a - 2b}
		\end{align*}
		so that we find the Euler classes $e_{35}$ and $e_{36}$ in the statement. 
		
		For the fixed point supported on $[1:0:0:0]$ and $[0:1:0:0]$ of type (3), filling in the tangent weights of $U_0$ in the last line of $V_{61}$ and switching $s_1$ and $s_3$ gives
		\begin{align*}
			t^{3a + 3b} &+ t^{-a-3b} + t^{2b} - t^{-3b+a} - t^{4a + 2b} - t^{5a +3b}.
		\end{align*} 
		This gives the Euler classes $e_{37}$ and $e_{38}$ from the statement.
		
		Again, one can check that for all possible congruence classes of $a,b$ modulo $4$, no signs are introduced because of the use of Proposition \ref{Marcs computation of euler classes}. 
	\end{proof} 
	The final class of ideal sheaves are those of which the ideal locally looks like the square of the ideal generated by two of the variables together with the remaining variable, i.e. they are supported on $[1:0:0:0]$ and $[0:1:0:0]$, and of one of the following types:
	\begin{enumerate} 
		\item $((x,y)^2,z)$ on $U_0$ and $((u,w)^2,v)$ on $U_1$. 
		\item $((y,z)^2,x)$ on $U_0$ and $((v,w)^2,u)$ on $U_1$. 
		\item $((x,z)^2,y)$ on $U_0$ and $((u,v)^2,w)$ on $U_1$. 
	\end{enumerate}
	Again, each of those have a corresponding fixed point supported on $[0:0:1:0]$ and $[0:0:0:1]$. 
	\begin{lem}
		The ideals of type (1) yield Euler classes 
		$$e_{39} = \frac{(3a+b)(3a-b)(5a-b)(2a-b)(2a+b)}{b^2(a-b)^2(a+b)}$$
		and	
		$$e_{310} = \frac{(3b+a)(3b-a)(5b-a)(2b-a)(2b+a)}{a^2(a-b)^2(a+b)}.$$
		The ideals of type (2) correspond to the Euler classes 
		$$e_{311} = \frac{9(3a+b)^3(3a-b)^3}{(a+b)^2(a-b)^2(a+3b)(a-3b)}$$
		and
		$$e_{312} = \frac{9(3b+a)^3(3b-a)^3}{(a+b)^2(a-b)^2(b+3a)(b-3a)}.$$
		The ideals of type (3) give rise to the Euler classes 
		$$e_{313} = \frac{(3a-b)(3a+b)(5a+b)(2a+b)(2a-b)}{b^2(a-b)(a+b)^2}$$
		and 
		$$e_{314} = \frac{(3b-a)(3b+a)(5b+a)(2b+a)(2b-a)}{a^2(b-a)(a+b)^2}.$$
	\end{lem}
	\begin{proof} 
		For $((x,y)^2,z)$ on $U_0$, plugging $Q = 1+s_1 + s_2$ into the formula from Proposition \ref{formula for trace of our representation} yields
		\begin{align*}
			V_{62} &= 2s_1^{-1} + 2s_2^{-1} + s_3^{-1} - s_1^{-1}s_2^{-1} - 2s_1^{-1}s_3^{-1} - 2s_2^{-1}s_3^{-1} \\
			&\quad + s_1s_3^{-1} + s_2s_3^{-1} - s_1^{-2}s_2^{-1} - s_1^{-1}s_2^{-2} \\
			&\quad + s_1^{-2}s_2 + s_1s_2^{-2} - s_1^{-2}s_2s_3^{-1} - s_1s_2^{-2}s_3^{-1}.
		\end{align*}
		For the fixed point of type (1), filling in the tangent weights of $U_0$ in $V_{62}$ gives 
		\begin{align*}
			2t^{2a} &+ 2t^{a-b} + t^{a+b} - 2t^{2a} - 2t^{3a+b} - t^{3a-b} \\
			+ t^{b-a} &+ t^{2b} - t^{5a-b} - t^{4a-2b} \\
			+ t^{3a+b} &+ t^{-2b} - t^{4a+2b} - t^{a-b}
		\end{align*}
		which yields the Euler classes $e_{39}$ and $e_{310}$ as in the statement. For the fixed point of type (2), plugging the tangent weights of $U_0$ into $V_{62}$ and switching $s_2$ and $s_3$ gives 
		\begin{align*}
			t^{a-b} &+ 2t^{2a} + 2t^{a+b} - 2t^{3a-b} - 2t^{2a} - t^{3a+b} \\
			+ t^{-a-b} &+ t^{-2b} - t^{5a+b} - t^{4a+2b} \\
			- t^{4a -2b} &- t^{a+b} + t^{3a-b} + t^{2b}
		\end{align*}
		which gives us the Euler classes $e_{311}$ and $e_{312}$ from the statement. Finally, for the fixed point of type (3), plugging the tangent weights of $U_0$ into $V_{62}$ and switching $s_1$ and $s_3$ gives 
		\begin{align*}
			t^{2a} &+ 2t^{a+b} + 2t^{a-b} - 2t^{3a+b} - 2t^{3a-b} - t^{2a} \\
			+ t^{a-b} &+ t^{a+b} - t^{3a+b} - t^{3a-b} \\
			-t^{3a+3b} &- t^{3a-3b} + t^{a+3b} + t^{a-3b} 
		\end{align*}
		which yields the Euler classes $e_{313}$ and $e_{314}$ as in the statement.
		
		One can check that for all possible congruence classes of $a,b$ modulo $4$, no signs are introduced because of the use of Proposition \ref{Marcs computation of euler classes}. 
	\end{proof} 
	Adding the Euler classes together proves Proposition \ref{computation Itilde6}. Combining Proposition \ref{computation Itilde6} with the Proposition \ref{computation n=2} and Proposition \ref{computation n=4}, we have proven Theorem \ref{main theorem}. 
	
	\phantomsection 
	\addcontentsline{toc}{section}{Bibliography}
	\bibliographystyle{plain}
	\bibliography{DTSources}
	
	\newpage 
	\noindent Anna M. Viergever \\
	Leibniz Universit\"at Hannover\\
	Fakult\"at f\"ur Mathematik, Welfengarten 1, 30167 Hannover, Germany \\
	E-Mail: \href{mailto:viergever@math.uni-hannover.de}{ viergever@math.uni-hannover.de} \\ \\ 
	Keywords: motivic homotopy theory, refined enumerative geometry, Donaldson-Thomas invariants, other fields\\Mathematics Subject Classification: 14N35, 14F42.
\end{document}